\documentclass[12pt]{article}
\usepackage{amssymb,latexsym,amsmath}

\begin{document}
\pagestyle{plain} \headheight=5mm  \topmargin=-5mm

\title{ Some relations between the topological and geometric filtration for smooth projective varieties   }
\author{ Wenchuan Hu }

\maketitle
\newtheorem{Def}{Definition}[section]
\newtheorem{Th}{Theorem}[section]
\newtheorem{Prop}{Proposition}[section]
\newtheorem{Not}{Notation}[section]
\newtheorem{Lemma}{Lemma}[section]
\newtheorem{Rem}{Remark}[section]
\newtheorem{Cor}{Corollary}[section]

\def\s{\section}
\def\ss{\subsection}

\def\d{\begin{Def}}
\def\t{\begin{Th}}
\def\p{\begin{Prop}}
\def\n{\begin{Not}}
\def\la{\begin{Lemma}}
\def\r{\begin{Rem}}
\def\c{\begin{Cor}}
\def\ee{\begin{equation}}
\def\aa{\begin{eqnarray}}
\def\ya{\begin{eqnarray*}}
\def\bd{\begin{description}}

\def\ed{\end{Def}}
\def\et{\end{Th}}
\def\epo{\end{Prop}}
\def\en{\end{Not}}
\def\el{\end{Lemma}}
\def\er{\end{Rem}}
\def\ec{\end{Cor}}
\def\eee{\end{equation}}
\def\eaa{\end{eqnarray}}
\def\ey{\end{eqnarray*}}
\def\ebd{\end{description}}

\def\nn{\nonumber}
\def\bp{{\bf Proof.}\hspace{2mm}}
\def\qe{\hfill$\Box$}
\def\lj{\langle}
\def\rj{\rangle}
\def\dd{\diamond}
\def\ox{\mbox{}}
\def\lb{\label}
\def\rel{\;{\rm rel.}\;}
\def\vp{\varepsilon}
\def\ep{\epsilon}
\def\mod{\;{\rm mod}\;}
\def\exp{{\rm exp}\;}
\def\Lie{{\rm Lie}}
\def\dim{{\rm dim}}
\def\im{{\rm im}\;}
\def\Lag{{\rm Lag}}
\def\Gr{{\rm Gr}}
\def\span{{\rm span}}
\def\Spin{{\rm Spin}}
\def\sign{{\rm sign}\;}
\def\Supp{{\rm Supp}\;}
\def\Sp{{\rm Sp}\;}
\def\ind{{\rm ind}\;}
\def\rank{{\rm rank}\;}
\def\Sg{{\Sp(2n,\C)}}
\def\Na{{\cal N}}
\def\det{{\rm det}\;}
\def\dist{{\rm dist}}
\def\deg{{\rm deg}}
\def\Tr{{\rm Tr}\;}
\def\ker{{\rm ker}\;}
\def\Vect{{\rm Vect}}
\def\H{{\bf H}}
\def\K{{\rm K}}
\def\R{{\bf R}}
\def\C{{\bf C}}
\def\Z{{\bf Z}}
\def\N{{\bf N}}
\def\F{{\bf F}}
\def\Da{{\bf D}}
\def\A{{\bf A}}
\def\La{{\bf L}}
\def\x{{\bf x}}
\def\y{{\bf y}}
\def\Ga{{\cal G}}
\def\Ha{{\cal H}}
\def\L{{\cal L}}
\def\Pa{{\cal P}}
\def\Ua{{\cal U}}
\def\E{{\rm E}}
\def\J{{\cal J}}

\def\m{{\rm m}}
\def\ch{{\rm ch}}
\def\gl{{\rm gl}}
\def\Gl{{\rm Gl}}
\def\Sp{{\rm Sp}}
\def\sf{{\rm sf}}
\def\U{{\rm U}}
\def\O{{\rm O}}
\def\F{{\rm F}}
\def\P{{\rm P}}
\def\D{{\rm D}}
\def\T{{\rm T}}
\def\Sa{{\rm S}}

\begin{center}{\bf \tableofcontents}\end {center}

\begin{center}{\bf Abstract}\end {center}{\hskip .2 in}
In the first part of this paper, we show that the assertion
``$T_pH_{k}(X,{\mathbb{Q}})=G_pH_k(X,{\mathbb{Q}})$" (which is
called the Friedlander-Mazur conjecture) is a birationally invariant
statement for smooth projective varieties $X$ when $p=\dim(X)-2$ and
when $p=1$. We also establish the  Friedlander-Mazur conjecture in
certain dimensions. More precisely, for a smooth projective variety
$X$, we show that the topological filtration $T_pH_{2p+1}(X,
{\mathbb{Q}})$ coincides with the geometric filtration
$G_pH_{2p+1}(X,{\mathbb{Q}})$ for all $p$. (Friedlander and Mazur
had previously shown that $T_pH_{2p}(X, {\mathbb{Q}})=G_pH_{2p}(X,
{\mathbb{Q}})$).   As a corollary, we conclude that for a smooth
projective threefold $X$,
$T_pH_k(X,{\mathbb{Q}})=G_pH_k(X,{\mathbb{Q}})$ for all $k\geq
2p\geq 0$ except for the case $p=1,k=4$. Finally, we show that the
topological and geometric filtrations always coincide if Suslin's
conjecture holds.

\s {Introduction}{\hskip .2 in} In this paper, all varieties are
defined over $\mathbb{C}$. Let $X$ be a  projective variety with
dimension $n$. Let ${\cal Z}_p(X)$ be the space of algebraic
$p$-cycles.

The  \textbf{Lawson homology} $L_pH_k(X)$ of $p$-cycles is defined
by $$L_pH_k(X) = \pi_{k-2p}({\cal Z}_p(X)) \quad for\quad k\geq
2p\geq 0,$$ where ${\cal Z}_p(X)$ is provided with a natural
topology (cf. \cite{Friedlander1}, \cite{Lawson1}). For general
background, the reader is referred to Lawson' survey paper
\cite{Lawson2}.

In \cite{Friedlander-Mazur}, Friedlander and Mazur showed that there
are  natural maps, called \textbf{cycle class maps}
 $$ \Phi_{p,k}:L_pH_{k}(X)\rightarrow H_{k}(X). $$

{\Def $$
\begin{array}{ll}
&L_pH_{k}(X)_{hom}:={\rm ker}\{\Phi_{p,k}:L_pH_{k}(X)\rightarrow
H_{k}(X)\};\\
&T_pH_{k}(X):={\rm Image}\{\Phi_{p,k}:L_pH_{k}(X)\rightarrow
H_{k}(X)\};\\
&T_pH_{k}(X,{\mathbb{Q}}):=T_pH_{k}(X)\otimes {\mathbb{Q}}.
\end{array}
 $$}

\medskip
It was shown in [\cite{Friedlander-Mazur}, \S 7] that the subspaces
$T_pH_k(X,{\mathbb{Q}})$ form a decreasing filtration:

$$\cdots\subseteq T_pH_k(X,{\mathbb{Q}})\subseteq T_{p-1}H_k(X,{\mathbb{Q}})
\subseteq\cdots\subseteq
T_0H_k(X,{\mathbb{Q}})=H_k(X,{\mathbb{Q}})$$ and
$T_pH_k(X,{\mathbb{Q}})$ vanishes if $2p>k$.

{\Def {\rm (\cite{Friedlander-Mazur})} Denote by
$G_pH_k(X,{\mathbb{Q}})\subseteq H_k(X,{\mathbb{Q}})$ the
$\mathbb{Q}$-vector subspace of $H_k(X,{\mathbb{Q}})$ generated by
the images of mappings $H_k(Y,{\mathbb{Q}})\rightarrow
H_k(X,{\mathbb{Q}})$, induced from all morphisms $Y\rightarrow X$ of
varieties of dimension $\leq k-p$.

The subspaces $G_pH_k(X,{\mathbb{Q}})$ also form a decreasing
filtration (called \textbf{geometric filtration}):
$$\cdots\subseteq G_pH_k(X,{\mathbb{Q}})\subseteq G_{p-1}H_k(X,{\mathbb{Q}})
\subset\cdots\subseteq G_0H_k(X,{\mathbb{Q}})\subseteq
H_k(X,{\mathbb{Q}})$$}

If $X$ is smooth, the Weak Lefschetz Theorem implies that
$G_0H_k(X,{\mathbb{Q}})= H_k(X,{\mathbb{Q}})$. Since $H_k(Y,{
\mathbb{Q}})$ vanishes for $k$ greater than twice the dimension of
$Y$, $G_pH_k(X,{\mathbb{Q}})$ vanishes if $2p>k$.

The following results have been proved by Friedlander and Mazur in
[FM]:

{\Th {\rm (\cite{Friedlander-Mazur})} Let $X$ be any projective
variety.
\begin{enumerate}
 \item  For non-negative integers $p$ and $k$,
$$T_pH_k(X,{\mathbb{Q}})\subseteq G_pH_k(X,{\mathbb{Q}}).$$
\item When $k= 2p$, $$T_pH_{2p}(X,{\mathbb{Q}})= G_pH_{2p}(X,{\mathbb{Q}}).$$
\end{enumerate}
}

\noindent\textbf{Question} (\cite{Friedlander-Mazur},
\cite{Lawson2}): Does one have equality in Theorem 1.1 when $X$ is a
smooth projective variety?

\medskip

Friedlander \cite{Friedlander2} has the following result: {\Th {\rm
(\cite{Friedlander2})} Let $X$ be a smooth projective variety of
dimension $n$. Assume that Grothendieck's Standard Conjecture B
(\cite{Grothendieck}) is valid for a resolution of singularities of
each irreducible subvariety of $Y\subset X$ of dimension $k-p$, then
$$T_pH_{k}(X,{\mathbb{Q}})=G_pH_k(X,{\mathbb{Q}}).$$ }

{\Rem{\rm (\cite{Lewis},\S 15.32)} The Grothendieck's Standard
Conjecture B is known to hold for a smooth projective variety $X$ in
the following cases:\rm
\begin{enumerate}
\item $\dim X\leq 2$.
\item Flag manifolds $X$.
\item Smooth complete intersections $X$.
\item Abelian varieties (due to D. Lieberman \cite{Lieberman}).

\end{enumerate}
}

\medskip
In this paper, we will use the tools in Lawson homology and the
methods given in [H] to show the following main results:

{\Th Let $X$ be a smooth projective variety of dimension $n$. If the
conclusion in Theorem 1.2 holds (without the assumption of
Grothendieck's Standard Conjecture B) for $X$ with
$p=1$,(resp.$p=n-2$) ($k$ arbitrary), then it also holds for
 any smooth projective variety $X^{\prime}$ which is
birationally equivalent to $X$ with $p=1$,(resp.$p=n-2$).}

{\Th For any smooth projective variety $X$, $$T_pH_{2p+1}(X,{
\mathbb{Q}})= G_pH_{2p+1}(X,{\mathbb{Q}}).$$ }

As corollaries, we have
 {\Cor Let $X$ be a smooth projective
3-fold. We have $T_pH_k(X,{\mathbb{Q}})=G_pH_k(X,{\mathbb{Q}})$ for
all
 $k\geq 2p\geq 0$ except for the case $p=1,k=4$. }

{\Cor Let $X$ be a smooth projective 3-fold with $H^{2,0}(X)=0$.
Then $T_pH_k(X,{\mathbb{Q}})=G_pH_k(X,{\mathbb{Q}})$ for any $k\geq
2p\geq 0$. In particular, it holds for $X$  a smooth hypersurface
and a complete intersection of dimension 3. }

\medskip
By using the K\"unneth formula in homology with rational
coefficient, we have

{\Cor Let $X$ be the product of  a smooth projective curve  and a
smooth simply connected projective surface. Then
$T_pH_k(X,{\mathbb{Q}})=G_pH_k(X,{\mathbb{Q}})$ for any $k\geq
2p\geq 0$. }

{\Cor For 4-folds $X$,  the assertion that
$T_pH_k(X,{\mathbb{Q}})=G_pH_k(X,{\mathbb{Q}})$ holds for all $k\geq
2p\geq 0$ is a birational invariant statement. In particular, if $X$
is a rational manifold with ${\rm dim}(X)\leq 4$, then the
conclusion in Theorem 1.2 holds for any $k\geq 2p \geq 0$ without
assumption of Grothendieck's Standard Conjecture B .

}


{\Rem A Conjecture given by Suslin (see
\cite{Friedlander-Haesemesyer-Walker}, \S 7) implies that
$L_pH_{n+p}(X^n)\cong H_{n+p}(X^n)$.}

\medskip
As an application of Theorem 1.4 and Proposition 3.1, we have the
following result:

{\Cor If the Suslin's Conjecture is true, then the topological
filtration is the same as the geometric filtration for a smooth
projective variety. }

\medskip
The main tools to prove this result are: the long exact localization
sequence given by Lima-Filho in \cite{Lima-Filho}, the explicit
formula for Lawson homology of codimension-one cycles on a smooth
projective manifold given by Friedlander in \cite{Friedlander1},
(and its generalization to general irreducible varieties, see
below), and the weak factorization theorem proved by Wlodarczyk in
\cite{Wlodarczyk} and in \cite{AKMW}.

\s {The Proof of the  Theorem 1.3} {\hskip .2 in} Let $X$ be a
smooth projective manifold of dimension $n$ and
$i_0:Y\hookrightarrow X$ be a smooth subvariety of codimension
$r\geq 2$. Let $\sigma:\tilde{X}_Y\rightarrow X$ be the blowup of
$X$ along $Y$, $\pi:D=\sigma^{-1}(Y)\rightarrow Y$ the nature map,
and $i:D=\sigma^{-1}(Y)\hookrightarrow \tilde{X}_Y$ the exceptional
divisor of the blowup. Set $U:= X-Y\cong \tilde{X}_Y - D$. Denote by
$j_0$ the inclusion $U\subset X$ and $j$ the inclusion $U\subset
\tilde{X}_Y$.

Now I list the Lemmas and Corollaries given in \cite{author}.

{\Lemma For each $p\geq 0$, we have the following commutative
diagram
$$
\begin{array}{ccccccccccc}
\cdots\rightarrow & L_pH_k(D) & \stackrel{i_*}{\rightarrow} &
L_pH_k({\tilde{X}_Y}) & \stackrel{j^*}{\rightarrow} & L_pH_k(U) &
\stackrel{\delta_*}{\rightarrow} & L_pH_{k-1}(D) & \rightarrow &
\cdots &
 \\

 & \downarrow \pi_*&   & \downarrow \sigma_*&   & \downarrow  \cong&   & \downarrow \pi_*&   &  &\\

 \cdots\rightarrow & L_pH_k(Y) & \stackrel{(i_0)*}{\rightarrow} &
L_pH_k({X}) & \stackrel{j_0^*}{\rightarrow}&
L_pH_k(U)&\stackrel{(\delta_0)_*}{\rightarrow} & L_pH_{k-1}(Y) &
\rightarrow & \cdots&

\end{array}
$$

}

{\Rem Since $\pi_* $ is surjective (there is an explicitly formula
for the Lawson homology of $D$, i.e., the Projective Bundle Theorem
proved by Friedlander and Gabber, see \cite{Friedlander-Gabber}), it
is easy to see that $\sigma_*$ is surjective.}

{\Cor If $p=0$, then we have the following commutative diagram
$$
\begin{array}{ccccccccccc}
\cdots\rightarrow & H_k(D) & \stackrel{i_*}{\rightarrow} &
H_k({\tilde{X}_Y}) & \stackrel{j^*}{\rightarrow} & H_k^{BM}(U) &
\stackrel{\delta_*}{\rightarrow} & H_{k-1}(D) & \rightarrow & \cdots
& \\

& \downarrow \pi_*&   & \downarrow \sigma_*&   & \downarrow  \cong&   & \downarrow \pi_*&   &  &\\

 \cdots\rightarrow & H_k(Y) & \stackrel{(i_0)*}{\rightarrow} &
H_k({X}) & \stackrel{j_0^*}{\rightarrow} &
H_k^{BM}(U)&\stackrel{(\delta_0)_*}{\rightarrow} & H_{k-1}(Y) &
\rightarrow & \cdots&

\end{array}
$$
Moreover, if $x\in H_k(D)$ maps to zero under $\pi_*$ and $i_*$,
then $x=0\in H_k(D)$.
 }

{\Cor If $p=n-2$, then we have the following commutative diagram
$$
\begin{array}{ccccccccccc}
\cdots\rightarrow & L_{n-2}H_k(D) & \stackrel{i_*}{\rightarrow} &
L_{n-2}H_k({\tilde{X}_Y}) & \stackrel{j^*}{\rightarrow} &
L_{n-2}H_k(U) & \stackrel{\delta_*}{\rightarrow} & L_{n-2}H_{k-1}(D)
& \rightarrow & \cdots & \\

& \downarrow \pi_*&   & \downarrow \sigma_*&   & \downarrow  \cong&
& \downarrow \pi_*&   &  &\\

 \cdots\rightarrow & L_{n-2}H_k(Y) & \stackrel{(i_0)*}{\rightarrow} &
L_{n-2}H_k({X}) & \stackrel{j_0^*}{\rightarrow} &
L_{n-2}H_k(U)&\stackrel{(\delta_0)_*}{\rightarrow} &
L_{n-2}H_{k-1}(Y) & \rightarrow & \cdots&

\end{array}
$$
}

{\Lemma For each $p\geq 0$, we have the following commutative
diagram
$$
\begin{array}{ccccccccccc}
\cdots\rightarrow & L_pH_k(D) & \stackrel{i_*}{\rightarrow} &
L_pH_k({\tilde{X}_Y}) & \stackrel{j^*}{\rightarrow} & L_pH_k(U) &
\stackrel{\delta_*}{\rightarrow} & L_pH_{k-1}(D) & \rightarrow &
\cdots & \\

& \downarrow \Phi_{p,k}&   & \downarrow \Phi_{p,k}&   & \downarrow
\Phi_{p,k}&   & \downarrow \Phi_{p,k-1}&   &  &\\

 \cdots\rightarrow & H_k(D) & \stackrel{i_*}{\rightarrow} &
H_k({\tilde{X}_Y}) & \stackrel{j^*}{\rightarrow} &
H_k^{BM}(U)&\stackrel{\delta_*}{\rightarrow} & H_{k-1}(D) &
\rightarrow & \cdots&

\end{array}
$$
In particular, it is true for $p=1, n-2$. }

\medskip
\bp See \cite{Lima-Filho} and also \cite{Friedlander-Mazur}.

\qe

{\Lemma For each $p\geq 0$, we have the following commutative
diagram
$$
\begin{array}{ccccccccccc}
\cdots\rightarrow & L_pH_k(Y) & \stackrel{(i_0)_*}{\rightarrow} &
L_pH_k({{X}}) & \stackrel{j^*}{\rightarrow}& L_pH_k(U) &
\stackrel{(\delta_0)_*}{\rightarrow} & L_pH_{k-1}(Y) & \rightarrow &
\cdots & \\

& \downarrow \Phi_{p,k}&   & \downarrow \Phi_{p,k}&   & \downarrow
\Phi_{p,k}&   & \downarrow \Phi_{p,k-1}&   &  &\\

 \cdots\rightarrow & H_k(Y) & \stackrel{(i_0)_*}{\rightarrow} &
H_k({X}) & \stackrel{j^*}{\rightarrow} &
H_k^{BM}(U)&\stackrel{(\delta_0)_*}{\rightarrow} & H_{k-1}(Y) &
\rightarrow & \cdots&

\end{array}
$$
In particular, it is true for $p=1, n-2$.}

\medskip
\bp See \cite{Lima-Filho} and also \cite{Friedlander-Mazur}. \qe

{\Rem The smoothness of $X$ and $Y$ is not necessary in the Lemma
2.3.}

{\Rem All the commutative diagrams of long exact sequences above
remain commutative and exact when tensored with ${\mathbb{Q}}$. We
will use these Lemmas and Corollaries with rational coefficients. }

\medskip
The following result will be used  several times in the proof of our
main theorem:

{\Th (Friedlander {\rm \cite{Friedlander1})} Let $W$ be any smooth
projective variety of dimension $n$. Then we have the following
isomorphisms
$$
\left\{
\begin{array}{l}
 L_{n-1}H_{2n}(W)\cong \Z,\\
 L_{n-1}H_{2n-1}(W)\cong H_{2n-1}(X,\Z),\\
 L_{n-1}H_{2n-2}(W)\cong H_{n-1,n-1}(X,\Z)=NS(W)\\
L_{n-1}H_{k}(X)=0 \quad for\quad k> 2n.\\

\end{array}
\right.
$$}

\noindent \textbf{The proof of Theorem 1.3 ( $p=n-2$ ):}

\medskip
There are two cases:
\medskip

\noindent\textbf{Case 1.} If $T_pH_k(X,{\mathbb{Q}})=G_pH_k(X,{
\mathbb{Q}})$, then
$T_pH_k({\tilde{X}_Y},{\mathbb{Q}})=G_pH_k({\tilde{X}_Y},{\mathbb{Q}})$.

\medskip
The injectivity of $T_pH_k({\tilde{X}_Y},{\mathbb{Q}})\rightarrow
G_pH_k({\tilde{X}_Y},{\mathbb{Q}})$ has been proved by Friedlander
and Mazur in \cite{Friedlander-Mazur}. We only need to show the
surjectivity.  Note that the case for $k=2p+1$ holds for any smooth
projective variety (Theorem 1.4). We only need to consider the cases
where $k\geq 2p+2$. In these cases, $k-p\geq p+2=n$, from the
definition of the geometric filtrations,  we have
$G_pH_k(\tilde{X},{\mathbb{Q}})=H_k({\tilde{X}_Y},{\mathbb{Q}})$ and
$G_pH_k(X,{\mathbb{Q}})=H_k(X,{\mathbb{Q}})$.

Let $b\in G_pH_k({\tilde{X}_Y},{\mathbb{Q}})$, and $a$ be the image
of $b$ under the the map  $\sigma_*:
H_k({\tilde{X}_Y},{\mathbb{Q}})\rightarrow H_k(X,{\mathbb{Q}})$,
i.e., $\sigma_*(b)=a$. By assumption, there exists an element
$\tilde{a}\in L_{n-2}H_k(X)\otimes{\mathbb{Q}}$ such that
$\Phi_{n-2,k}(\tilde{a})=a$. Since $\sigma_*:
L_{n-2}H_k({\tilde{X}_Y})\otimes{\mathbb{Q}}\rightarrow
L_{n-2}H_k(X)\otimes{\mathbb{Q}}$ is surjective (\cite{author}),
there exists an element $\tilde{b}\in
L_{n-2}H_k(X)\otimes{\mathbb{Q}}$ such that
$\sigma_*(\tilde{b})=\tilde{a}$. By the following commutative
diagram

$$
\begin{array}{cccccc}
\    &   & \  & L_{n-2}H_k({\tilde{X}_Y})\otimes{\mathbb{Q}}
& \stackrel{\sigma_*}{\rightarrow} & L_{n-2}H_k(X)\otimes{\mathbb{Q}}  \\

&   &   & \downarrow \Phi_{n-2,k}&   & \downarrow  \Phi_{n-2,k}   \\

&  &   & H_k({\tilde{X}_Y},{\mathbb{Q}}) &
\stackrel{\sigma_*}{\rightarrow} & H_k(X,{\mathbb{Q}}),

\end{array}
$$
we have $\Phi_{n-2,k}(\tilde{b})-b$ maps to zero in
$H_k(X,{\mathbb{Q}})$. By the commutative diagram in Corollary 2.1,
$j^*(\Phi_{n-2,k}(\tilde{b})-b)=0\in H_k^{BM}(U,{\mathbb{Q}})$. From
the exactness of the upper long exact sequence in Corollary 2.1,
there exists an element $c\in H_k(D,{\mathbb{Q}})$ such that
$i_*(c)=\Phi_{n-2,k}(\tilde{b})-b$. From Theorem 2.1, we find that
$\Phi_{n-2,k}:L_{n-2}H_{k}(D)\otimes{\mathbb{Q}}\rightarrow
H_{k}(D)\otimes{\mathbb{Q}}$ is an isomorphism for $k\geq 2n-2$.
Hence there exists an element $\tilde{c}\in
L_{n-2}H_{k}(D)\otimes{\mathbb{Q}}$ such that
$i_*(\Phi_{n-2,k}(\tilde{c}))=\Phi_{n-2,k}(\tilde{b})-b$. Therefore
$\Phi_{n-2,k}(\tilde{b}-i_*(\tilde{c}))=b$, i.e., the surjectivity
of $T_pH_k({\tilde{X}_Y},{\mathbb{Q}})\rightarrow
G_pH_k({\tilde{X}_Y},{\mathbb{Q}})$.

\medskip
On the other hand, we need to show

\medskip
\noindent \textbf{Case 2.} If
$T_pH_k(\tilde{X}_Y,{\mathbb{Q}})=G_pH_k(\tilde{X}_Y,{\mathbb{Q}})$,
then $T_pH_k(X,{\mathbb{Q}})=G_pH_k(X,{\mathbb{Q}})$.

\medskip
This part is relatively easy. By Theorem 1.4, we only need to
consider the cases that $k\geq 2p+2=2n-2$. Let $a\in G_pH_k(X,{
\mathbb{Q}})=H_k(X,{\mathbb{Q}})$. From the blow up formula for
singular homology (cf. [GH]), we know
$\sigma_*:H_k(\tilde{X}_Y,{\mathbb{Q}})\rightarrow
H_k(X,{\mathbb{Q}})$ is surjective. Then there exists an element
$b\in H_k(\tilde{X}_Y,{\mathbb{Q}})$ such that $\sigma_*(b)=a$. By
assumption, we can find an element $\tilde{b}\in
L_{n-2}H_k(\tilde{X}_Y,{\mathbb{Q}})$ such that
$\Phi_{n-2,k}(\tilde{b})=b$. Set $\tilde{a}=\sigma_*(\tilde{b})$.
Then $\Phi_{n-2,k}(\tilde{a})=a$ under the natural map
$\Phi_{n-2,k}$. This is exactly the surjectivity we want.

This completes the proof for a blow-up along a smooth codimension at
least two subvariety $Y$ in $X$.

\qe

\noindent\textbf{The proof of Theorem 1.3 ($p=1$):}

\medskip
The injectivity of the map $T_1H_k(W,{\mathbb{Q}})\rightarrow
G_1H_k(W,{\mathbb{Q}})$ has been proved for any smooth projective
variety $W$ by Friedlander and Mazur in \cite{Friedlander-Mazur}. We
only need to show the surjectivity under certain assumption.

Similar to the case $p=n-2$, we also have two cases:

\medskip
\noindent\textbf{Case A.} If
$T_1H_k(X,{\mathbb{Q}})=G_1H_k(X,{\mathbb{Q}})$, then
$T_1H_k(\tilde{X}_Y,{\mathbb{Q}})=G_1H_k(\tilde{X}_Y,{\mathbb{Q}})$.

\medskip
From Theorem 1.4, the case where $k=3$ holds for any smooth
projective variety. We only need to consider the cases where $k\geq
4$.

Let $b\in G_1H_k(\tilde{X}_Y,{\mathbb{Q}})$. Denote by $a$ the image
of $b$ under the the map  $\sigma_*: H_k(\tilde{X}_Y,{\bf
Q})\rightarrow H_k(X,{\mathbb{Q}})$, i.e., $\sigma_*(b)=a$. From the
blow up formula for singular homology and the definition of the
geometric filtration, we have
$\sigma_*(G_1H_k(\tilde{X}_Y,{\mathbb{Q}}))=
G_1H_k(X,{\mathbb{Q}})$.




By assumption, there exists an element $\tilde{a}\in
L_{1}H_k(X)\otimes{\mathbb{Q}}$ such that $\Phi_{1,k}(\tilde{a})=a$.
Since $\sigma_*: L_{1}H_k(\tilde{X}_Y)\otimes{\mathbb{Q}}\rightarrow
L_{1}H_k(X)\otimes{\mathbb{Q}}$ is surjective ([H]), there exists an
element $\tilde{b}\in L_{1}H_k(\tilde{X}_Y)\otimes{\mathbb{Q}}$ such
that $\sigma_*(\tilde{b})=\tilde{a}$. By the following commutative
diagram
$$
\begin{array}{ccccccccccc}
\    &   & \  & L_{1}H_k({\tilde{X}_Y})\otimes{\mathbb{Q}}
 & \stackrel{\sigma_*}{\rightarrow} & L_{1}H_k(X)\otimes{\mathbb{Q}} &   &   &   &   & \\

 &   &   & \downarrow \Phi_{1,k}&   & \downarrow  \Phi_{1,k}&   &  &   &  &\\

     &  &   &
  H_k({\tilde{X}_Y},{\mathbb{Q}}) & \stackrel{\sigma_*}{\rightarrow} & H_k(X,{\mathbb{Q}}),&  &
  &   &  &

\end{array}
$$
we have $\Phi_{1,k}(\tilde{b})-b$ maps to zero in
$H_k(X,{\mathbb{Q}})$. By the commutative diagram in Corollary 2.1,
$j^*(\Phi_{1,k}(\tilde{b})-b)=0\in H_k^{BM}(U,{\mathbb{Q}})$. From
the exactness of the upper long exact sequence in Corollary 2.1,
there exists an element $c\in H_k(D,{\mathbb{Q}})$ such that
$i_*(c)=\Phi_{1,k}(\tilde{b})-b$. Set $d=\pi_*(c)\in
H_{k}(Y,{\mathbb{Q}})$. By the commutative diagram in Corollary 2.1,
$d$ maps to zero under $(i_0)_*: H_k(Y,{\mathbb{Q}})\rightarrow
H_k(X,{\mathbb{Q}})$. Hence there exists an element $e\in
H_{k+1}^{BM}(U,{\mathbb{Q}})$ such that whose image is $d$ under the
boundary map $(\delta_0)_*$. Let $\tilde{d}\in H_k(D,{\mathbb{Q}})$
be the image of $e$ under this boundary map $\delta_*:
H_{k+1}^{BM}(U,{\mathbb{Q}})\rightarrow H_k(D,{\mathbb{Q}})$.
Therefore, the image of $c-\tilde{d}$ is zero under $\pi_*$ in
$H_k(Y,{\mathbb{Q}})$  and  is also zero under $i_*$ in
$H_k(\tilde{X}_Y,{\mathbb{Q}})$. Note that $D$ is a bundle over $Y$
with projective spaces as fibers. From the ``projective bundle
theorem" for the singular homology (cf.[GH]), we have
$H_{k}(D,{\mathbb{Q}})\cong H_{k}(Y,{\mathbb{Q}})\oplus
H_{k-2}(Y,{\mathbb{Q}})\oplus \cdots\oplus
H_{k-2r+2}(Y,{\mathbb{Q}})$. From this,  we have $c-\tilde{d}\in
H_{k-2}(Y,{\mathbb{Q}})\oplus\cdots\oplus
H_{k-2r+2}(Y,{\mathbb{Q}})$. By the revised Projective Bundle
Theorem (\cite{Friedlander-Gabber}, and \cite{author} the revised
case essentially due to Complex Suspension Theorem \cite{Lawson1})
and Dold-Thom Theorem \cite{Dold-Thom}, we have
$L_1H_k(D,{\mathbb{Q}})\cong L_1H_k(Y,{\mathbb{Q}})\oplus
L_0H_{k-2}(Y,{\mathbb{Q}})\oplus\cdots\oplus
L_{2-r}H_{k-2r+2}(Y,{\mathbb{Q}})\cong L_1H_k(Y,{\mathbb{Q}})\oplus
H_{k-2}(Y,{\mathbb{Q}})\oplus\cdots\oplus
H_{k-2r+2}(Y,{\mathbb{Q}})$, where $r$ is the codimension of $Y$.
Since $c-\tilde{d}\in H_{k-2}(Y,{\mathbb{Q}})\oplus\cdots\oplus
H_{k-2r+2}(Y,{\mathbb{Q}})$ and
$L_0H_{k-2}(Y,{\mathbb{Q}})\oplus\cdots\oplus
L_{2-r}H_{k-2r+2}(Y,{\mathbb{Q}})\cong
H_{k-2}(Y,{\mathbb{Q}})\oplus\cdots\oplus
H_{k-2r+2}(Y,{\mathbb{Q}})$, there exists an element $f\in
L_1H_k(D,{\mathbb{Q}})$ such that $\Phi_{1,k}(f)=c-\tilde{d}$.
Therefore we obtain $\Phi_{1,k}(\tilde{b}-i_*(f))=b$. This is the
surjectivity we need.

\medskip
\noindent\textbf{Case B.} If
$T_1H_k(\tilde{X}_Y,{\mathbb{Q}})=G_1H_k(\tilde{X}_Y,{\mathbb{Q}})$,
then $T_1H_k(X,{\mathbb{Q}})=G_1H_k(X,{\mathbb{Q}})$.

\medskip
This part is also relatively easy. Note that $k\geq 4$. Let $a\in
G_1H_k(X,{\mathbb{Q}})\subseteq H_k(X,{\mathbb{Q}})$, then there
exists an element $b\in G_1H_k(\tilde{X}_Y,{\mathbb{Q}})$ such that
$\sigma_*(b)=a$. By assumption, we can find an element $\tilde{b}\in
L_{1}H_k(\tilde{X}_Y,{\mathbb{Q}})$ such that
$\Phi_{1,k}(\tilde{b})=b$. Set $\tilde{a}=\sigma_*(\tilde{b})$. Then
$\Phi_{1,k}(\tilde{a})=a$ under the natural transformation
$\Phi_{1,k}$. This is exactly the surjectivity in these cases.

This completes the proof for one blow-up along a smooth codimension
at least two subvariety $Y$ in $X$.

\qe

Now recall the weak factorization Theorem proved in \cite{AKMW} (and
also \cite{Wlodarczyk}) as follows:

{\Th (\cite{AKMW} Theorem 0.1.1, \cite{Wlodarczyk})
 Let $\varphi \colon X \to X'$ be a birational map of smooth
complete varieties over an algebraically closed field of
characteristic zero, which is an isomorphism over an open set $U$.
Then $f$ can be factored as a sequence of birational maps
$$X = X_0
\stackrel{\varphi_1}{\rightarrow} X_1
\stackrel{\varphi_2}{\rightarrow}\dots
\stackrel{\varphi_{n+1}}{\rightarrow} X_n = X'$$ where each $X_i$ is
a smooth complete variety, and $\varphi_{i+1}: X_i \to X_{i+1}$ is
either a blowing-up or a blowing-down of a smooth subvariety
disjoint from $U$. }

\qe

{\Rem From the proof of the Theorem 1.3, we can draw the following
conclusions:

\begin{enumerate}
\item If $$T_rH_{k}(Y,{\mathbb{Q}})=G_rH_k(Y,{\mathbb{Q}})$$ for all $k$ is
true for algebraic r-cycles with  $r\geq p$ for $\dim(Y)=n$, then
$$``T_{p-1}H_{k}(X,{\mathbb{Q}})=G_{p-1}H_k(X,{\mathbb{Q}}),\quad\forall k"$$ is
a birationally invariant statement for smooth projective varieties
$X$ with $\dim(X)\leq n+2$.

\item If $$T_rH_{k}(Y,{\mathbb{Q}})=G_rH_k(Y,{\mathbb{Q}})$$ for all $k$ is
true for r-algebraic cycles with  $r\leq p$ for $\dim(Y)=n$, then
$$``T_{p+1}H_{k}(X,{\mathbb{Q}})=G_{p+1}H_k(X,{\mathbb{Q}}),\quad \forall k"$$ is
a birationally invariant statement for smooth projective varieties
$X$ with  $\dim(X)\leq n+2$.

\end{enumerate}
}


\s {The Proof of the Theorem 1.4}

{\Prop For any irreducible projective variety $Y$ of dimension
$n$, we have
$$
\left\{
\begin{array}{l}
 L_{n-1}H_{2n}(X)\cong \Z,\\
 L_{n-1}H_{2n-1}(X)\cong H_{2n-1}(X,\Z),\\
 L_{n-1}H_{2n-2}(X)\rightarrow H_{2n-2}(X,\Z) {\quad is\quad injective},\\
L_{n-1}H_{k}(X)=0 \quad for\quad k> 2n.\\

\end{array}
\right.
$$}

{\Rem When $Y$ is smooth projective, Friedlander have drawn a
stronger conclusion, i.e., besides those in the proposition,
$L_{n-1}H_{2n-2}(Y)\cong H_{n-1,n-1}(X, \Z)=NS(X)$.}

\medskip
\bp Set $S=Sing(Y)$, the set of singular points. Then $S$ is the
union of proper irreducible subvarieties. Set $S= (\cup_{i}
S_i)\bigcup S'$, where $\dim (S_i)=n-1$ and $S'$ is the union of
subvarieties with dimension $\leq n-2$. Let $V=Y-S$ be the smooth
open part of $Y$. According to Hironaka \cite{Hironaka}, we can find
$\tilde{Y}$ such that $\tilde{Y}$ is a smooth compactification of
$V$. Let $D=\tilde{Y}-V$. $D$ is a divisor on $\tilde{Y}$ with
normal crossing. Denote by $i_0:S \hookrightarrow Y$ and $i:D
\hookrightarrow \tilde{Y}$ the inclusions of closed sets. Denote by
$j_0:V \hookrightarrow Y$ and  $j:V \hookrightarrow \tilde{Y}$ the
inclusions of open sets.

\medskip
There are a few cases:

\medskip
\noindent\textbf{ Case 1:} $k\geq 2n$.

\medskip
By the localization long exact sequence in Lawson homology
$$\cdots\rightarrow L_{n-1}H_k(S)\rightarrow L_{n-1}H_k(Y)\rightarrow L_{n-1}H_k(V)\rightarrow
L_{n-1}H_{k-1}(S)\rightarrow\cdots,$$
we have
$$L_{n-1}H_k(Y)\cong L_{n-1}H_k(V) \quad for\quad k\geq 2n$$
since $L_{n-1}H_k(S)=0$ for $k\geq 2n-1$.

By the localization exact sequence in homology

$$\cdots\rightarrow H_k(S)\rightarrow H_k(Y)\rightarrow H_k^{BM}(V)\rightarrow
H_{k-1}(S)\rightarrow\cdots,$$
we have
$$H_k(Y)\cong H_k^{BM}(V) \quad for\quad k\geq 2n$$
since $H_k(S)=0$ for $k\geq 2n-1$. Here $H_k^{BM}(V)$ is the
Borel-Moore homology.

Similarly, $$L_{n-1}H_k(\tilde{Y})\cong L_{n-1}H_k(V) \quad
for\quad k\geq 2n$$ and $$H_k(\tilde{Y})\cong H_k^{BM}(V) \quad
for\quad k\geq 2n.$$

Since $\tilde{Y}$ is smooth, we have $L_{n-1}H_k(\tilde{Y})\cong
H_k(\tilde{Y})$ for $k\geq 2n$(cf. \cite{Friedlander1}). This
completes the proof for the case $k\geq 2n$.

\medskip
\noindent \textbf{Case 2:}  $k=2n-1$.

\medskip
Applying Lemma 2.3 to the pair $(Y,S)$ for $p=n-1$, we have the
commutative diagram of the long exact sequence

\begin{eqnarray}
\begin{array}{ccccccccccc}
\ 0\rightarrow & L_{n-1}H_{2n-1}(Y) & \stackrel{j_0^*}{\rightarrow}
& L_{n-1}H_{2n-1}(V) & \stackrel{(\delta_0)*}{\rightarrow} &
L_{n-1}H_{{2n-2}}(S) & \stackrel{(i_0)*}{\rightarrow} &
L_{n-1}H_{2n-2}(Y) &
\rightarrow & \cdots & \\
& \downarrow \Phi_{n-1,2n-1} & & \downarrow \Phi_{n-1,2n-1}&
& \downarrow \Phi_{n-1,2n-2}& & \downarrow \Phi_{n-1,2n-2}& & &\\

\ 0\rightarrow & H_{2n-1}(Y) & \stackrel{j_0^*}{\rightarrow} &
H_{2n-1}^{BM}(V) & \stackrel{(\delta_0)*}{\rightarrow} &
H_{{2n-2}}(S)&\stackrel{(i_0)*}{\rightarrow} & H_{2n-2}(Y) &
\rightarrow & \cdots&

\end{array}
\end{eqnarray}

Similarly, applying Lemma 2.3 to the pair $(\tilde{Y},D)$ for
$p=n-1$, we have the commutative diagram of the long exact sequence

\begin{eqnarray}
\begin{array}{ccccccccccc}
\ 0\rightarrow & L_{n-1}H_{2n-1}(\tilde{Y}) &
\stackrel{j^*}{\rightarrow} & L_{n-1}H_{2n-1}(V) &
\stackrel{\delta_*}{\rightarrow}& L_{n-1}H_{{2n-2}}(D) &
\stackrel{i_*}{\rightarrow} & L_{n-1}H_{2n-2}(\tilde{Y}) &
\rightarrow & \cdots & \\
  & \downarrow \tilde{\Phi}_{n-1,2n-1}&& \downarrow \Phi_{n-1,2n-1}&
  & \downarrow\Phi_{n-1,2n-2}& & \downarrow\tilde{\Phi}_{n-1,2n-2}&   &  &\\

\ 0\rightarrow & H_{2n-1}(\tilde{Y}) & \stackrel{j^*}{\rightarrow} &
H_{2n-1}^{BM}(V) & \stackrel{\delta_*}{\rightarrow} &
H_{{2n-2}}(D)&\stackrel{i_*}{\rightarrow} & H_{2n-2}(\tilde{Y}) &
\rightarrow & \cdots&
\end{array}
\end{eqnarray}

Note that
$\tilde{\Phi}_{n-1,2n-2}:L_{n-1}H_{2n-2}(\tilde{Y})\rightarrow
H_{2n-2}(\tilde{Y})$ is injective,
$\tilde{\Phi}_{n-1,2n-1}:L_{n-1}H_{2n-1}(\tilde{Y})\cong
H_{2n-1}(\tilde{Y})$  and
$\tilde{\Phi}_{n-1,2n-2}:L_{n-1}H_{2n-2}(D)\cong H_{2n-2}(D)\cong
\Z^m$, where $m$ is the number of irreducible varieties of $D$. From
(2) and the Five Lemma, we have the isomorphism

\begin{eqnarray} {\Phi}_{n-1,2n-1}:L_{n-1}H_{2n-1}(V)\cong
H_{2n-1}^{BM}(V).
\end{eqnarray}

From (1), (3) and the Five Lemma, we have the following isomorphism
$${\Phi}_{n-1,2n-1}:L_{n-1}H_{2n-2}(Y)\cong H_{2n-2}(Y).$$

\medskip
\noindent\textbf{Case 3:}  $k=2n-2$.

\medskip
Now the commutative diagram (1) is rewritten in the following way:
\begin{eqnarray}
\begin{array}{ccccccccccc}
\ \cdots & \rightarrow & L_{n-1}H_{2n-1}(V) &
\stackrel{(\delta_0)_*}{\rightarrow}& L_{n-1}H_{{2n-2}}(S) &
\stackrel{(i_0)_*}{\rightarrow} & L_{n-1}H_{2n-2}(Y) &
\stackrel{j_0^*}{\rightarrow} & L_{n-1}H_{2n-2}(V) & \rightarrow&0\\
& & \downarrow \Phi_{n-1,2n-1}& & \downarrow \Phi_{n-1,2n-2}&
& \downarrow \Phi_{n-1,2n-2}& &\downarrow \Phi_{n-1,2n-2} & &\\

\ \cdots & \rightarrow & H_{2n-1}^{BM}(V) &
\stackrel{(\delta_0)_*}{\rightarrow} &
H_{{2n-2}}(S)&\stackrel{(i_0)_*}{\rightarrow} & H_{2n-2}(Y) &
\stackrel{j_0^*}{\rightarrow} & H_{2n-2}^{BM}(V)&\rightarrow &0

\end{array}
\end{eqnarray}

In the  commutative diagram (2), we can show that the injective maps

\begin{eqnarray}
 j^*: H_{2n-1}(\tilde{Y})\rightarrow H_{2n-1}^{BM}(V)
\end{eqnarray}
and
\begin{eqnarray}
 j^*: L_{n-1}H_{2n-1}(\tilde{Y})\rightarrow L_{n-1}H_{2n-1}(V)
\end{eqnarray}
are actually isomorphisms. Hence the commutative diagram (2) reduces
to the following diagram:

\begin{eqnarray}
\begin{array}{ccccccccccc}
\    & 0 & \rightarrow & L_{n-1}H_{{2n-2}}(D) & \rightarrow &
L_{n-1}H_{2n-2}(\tilde{Y}) & \rightarrow
 & L_{n-1}H_{{2n-2}}(V) & \rightarrow & 0 &\\
&  &  & \downarrow\Phi_{n-1,2n-2}& &
\downarrow\tilde{\Phi}_{n-1,2n-2}&   &\downarrow\Phi_{n-1,2n-2}  & &
&\\

\  & 0 & \rightarrow & H_{{2n-2}}(D)&\rightarrow &
H_{2n-2}(\tilde{Y}) & \rightarrow & H_{{2n-2}}^{BM}(V)&
\rightarrow & 0&
\end{array}
\end{eqnarray}
To see (5) are surjective, by the exactness of the rows in (2) we
only need to show that the maps $i_*: H_{{2n-2}}(D) \rightarrow
H_{2n-2}(\tilde{Y})$ are injective. Note that $\tilde{Y}$ is a
compact K\"ahlar manifold, and  the homology class of an algebraic
subvariety is nontrivial in the homology of the K\"ahlar manifold.
From these, we get the injectivity of $i_*$. The surjectivity of (6)
follows from the same reason.

\medskip
We need the following lemma.

{\Lemma The natural transformation
$\Phi_{n-1,2n-2}:L_{n-1}H_{2n-2}(V)\rightarrow H_{2n-2}^{BM}(V)$ is
injective.}

\medskip
\bp $a\in L_{n-1}H_{2n-2}(V)$ such that $\Phi_{n-1,2n-2}(a)=0 \in
H_{2n-2}^{BM}(V)$. Since the map
$j^*:L_{n-1}H_{2n-2}(\tilde{Y})\rightarrow L_{n-1}H_{2n-2}(V)$ is
surjective, there exists an element $b\in
L_{n-1}H_{2n-2}(\tilde{Y})$ such that $j^*(b)=a$. Set
$\tilde{b}=\Phi_{n-1,2n-2}(b)\in H_{2n-2}(\tilde{Y})$. By the
commutativity of the diagram, we have $j^*(\tilde{b})=0$ under the
map $j^*:H_{2n-2}(\tilde{Y})\rightarrow H_{2n-2}^{BM}(V)$. By the
exactness of the bottom row in the commutative diagram (7), there
exists an element $ \tilde{c}\in H_{{2n-2}}(D)$ such that the image
of $\tilde{c}$ under the map  $i_*: H_{{2n-2}}(D) \rightarrow
H_{2n-2}(\tilde{Y})$ is $\tilde{b}$. Now note that $\Phi_{n-1,2n-2}:
L_{n-1}H_{{2n-2}}(D)\rightarrow H_{{2n-2}}(D)$ is an isomorphism,
there exists an element $c\in L_{n-1}H_{{2n-2}}(D)$ such that
$\Phi_{n-1,2n-2}(c)=\tilde{c}$. Hence $\Phi_{n-1,2n-2}(i_*(c)-b)=0$.
Note that $\Phi_{n-1,2n-2}: L_{n-1}H_{{2n-2}}(\tilde{Y})\rightarrow
H_{{2n-2}}(\tilde{Y})$ is injective since $\tilde{Y}$ is smooth and
of dimension $n$ (cf. [F1]). Hence we get $i_*(c)=b$, i.e., $b$ is
in the image of the map $i_*:L_{n-1}H_{{2n-2}}(D)\rightarrow
L_{n-1}H_{{2n-2}}(\tilde{Y})$. Therefore $a=0$ by the exactness of
the top row of the commutative diagram (7).

\qe

We need to show that $\Phi_{n-1,2n-2}:L_{n-1}H_{2n-2}(Y)\rightarrow
H_{2n-2}(Y)$ is injective. For $a\in L_{n-1}H_{2n-2}(Y)$ such that
$\Phi_{n-1,2n-2}(a)=0 \in  H_{2n-2}(Y)$. By the commutative diagram
(4) and the Lemma 3.1,  the image of $a$ under
$j_0^*:L_{n-1}H_{2n-2}(Y)\rightarrow L_{n-1}H_{2n-2}(V)$ is zero.
Hence there exists an element $b\in L_{n-1}H_{2n-2}(S)$ such that
the image of $(i_0)_*:L_{n-1}H_{2n-2}(S)\rightarrow
L_{n-1}H_{2n-2}(Y)$ is $a$, i.e., $(i_0)_*(b)=a$. Set
$\tilde{b}=\Phi_{n-1,2n-2}(b)$. Then the image of $\tilde{b}$ under
the map $(i_0)_*: H_{2n-2}(S)\rightarrow H_{2n-2}(Y)$ is zero. By
exactness of the bottom row in the commutative diagram (4), there
exists an element $\tilde{c}$ such that its image under the map
$H_{2n-1}^{BM}(V)\rightarrow H_{2n-2}(S)$ is $\tilde{b}$. By the
result in \textbf{Case 2}, $\Phi_{n-1,2n-1}:
L_{n-1}H_{2n-1}(V)\rightarrow H_{2n-1}^{BM}(V)$ is an isomorphism.
Hence there exists an element $c\in L_{n-1}H_{2n-1}(V)$ such that
$\Phi_{n-1,2n-1}(c)=\tilde{c}$. Now since $\Phi_{n-1,2n-2}:
L_{n-1}H_{2n-2}(S)\rightarrow H_{2n-2}(S)$ is an isomorphism, the
image of $c$ under the map $L_{n-1}H_{2n-1}(V)\rightarrow
L_{n-1}H_{2n-2}(S)$ is exactly $b$. Now the exactness of the top row
of the commutative diagram (4) implies the vanishing of $a$.

The proof of the proposition is done.

 \qe

By using this proposition, we will give a proof of Theorem 1.4.

\medskip

\noindent \textbf{Proof of Theorem 1.4:}

\medskip
For any smooth projective variety $X$, the injectivity of
$T_pH_{2p+1}(X,{ \mathbb{Q}})\rightarrow
G_pH_{2p+1}(X,{\mathbb{Q}})$ has been proved in
[\cite{Friedlander-Mazur}, \S 7]. We only need to show the
surjectivity of $T_pH_{2p+1}(X,{ \mathbb{Q}})\rightarrow
G_pH_{2p+1}(X,{\mathbb{Q}})$.
 For any subvariety $i: Y\subset X$, we
denote by $V=:X-Y$ the complementary of $Y$ in $X$. We have the
following commutative diagram of the long exact sequences (Lemma
2.3, or \cite{Lima-Filho}):
$$
\begin{array}{ccccccccccc}
\cdots\rightarrow & L_pH_{2p+1}(Y) & \rightarrow & L_pH_{2p+1}(X)
& \rightarrow & L_pH_{2p+1}(V) & \rightarrow & L_pH_{2p}(Y) &
\rightarrow & \cdots & \\
& \downarrow \Phi_{p,2p+1}& & \downarrow \Phi_{p,2p+1}&
& \downarrow \Phi_{p,2p+1}& & \downarrow \Phi_{p,2p}&   &  &\\

\cdots\rightarrow &H_{2p+1}(Y) & \rightarrow & H_{2p+1}(X) &
\rightarrow & H_{2p+1}^{BM}(V)&\rightarrow & H_{2p}(Y) &
\rightarrow & \cdots&

\end{array}
$$
Obviously, the above commutative diagram holds when tensored  with
${\mathbb{Q}}$. In the following, we only consider the commutative
diagrams with  ${\mathbb{Q}}$-coefficient.

Now let $a\in G_pH_{2p+1}(X,{\mathbb{Q}})$, by definition, we can
assume that $a$ lies in the image of the map $i_*:
H_{2p+1}(Y,{\mathbb{Q}}) \rightarrow H_{2p+1}(X,{\mathbb{Q}})$ for
some subvariety $Y\subset X$ with dimension $\dim Y=(2p+1)-p=p+1$.
Hence there exists an element $b\in H_{2p+1}(Y,{\mathbb{Q}})$ such
that $i_*(b)=a$. By the Proposition 3.1, we know that
$\Phi_{p,2p+1}: L_pH_{2p+1}(Y)\otimes{\mathbb{Q}}\rightarrow
H_{2p+1}(Y,{\mathbb{Q}})$ is an isomorphism. Therefore there exists
an element $\tilde{b}\in L_pH_{2p+1}(Y)\otimes{\mathbb{Q}}$ such
that $\Phi_{p,2p+1}(\tilde{b})=b$. Set $\tilde{a}=i_*(\tilde{b})$.
Then $\tilde{a}$ maps to $a$ under the map
$L_pH_{2p+1}(X)\otimes{\mathbb{Q}}\rightarrow
H_{2p+1}(X,{\mathbb{Q}})$. By the definition of the topological
filtration, $a\in T_pH_{2p+1}(X,{\mathbb{Q}})$. This completes the
proof of surjectivity of $T_pH_{2p+1}(X,{ \mathbb{Q}})\rightarrow
G_pH_{2p+1}(X,{\mathbb{Q}})$.

 \qe

{\Rem In the proof of the surjectivity of Theorem 1.4, the
assumption of smoothness is not necessary, more precisely, for any
irreducible projective variety $X$, the image of the natural
transformation $\Phi_{p,2p+1}:
L_pH_{2p+1}(X,{\mathbb{Q}})\rightarrow H_{2p+1}(X,{\mathbb{Q}})$
contains $G_pH_{2p+1}(X,{\mathbb{Q}})$. }

{\Rem Independently, M. Warker  has recently also obtained this
result (\cite{Walker}, Prop. 2.5]).}

\medskip
Now we prove the corollaries 1.2-1.5.

\medskip
\noindent\textbf{The proof of Corollary 1.1:} By Theorem 1.1 and
1.4, Dold-Thom Theorem and Proposition 3.1, we only need to show the
cases that $p=1,k\geq 5$. Now the following commutative diagram
(\cite{Friedlander-Mazur}, Prop.6.3)
$$
\begin{array}{ccccc}
 & L_2H_k(X)\otimes {\mathbb{Q}}  &  \stackrel{s}{\rightarrow} & L_1H_k(X)\otimes {\mathbb{Q}}  &    \\

 &  \downarrow \Phi_{2,k} &   &\downarrow \Phi_{1,k} &\\

 & H_k(X,{\mathbb{Q}}) &\cong   & H_k(X,{\mathbb{Q}}). &
\end{array}
$$
shows that if $L_2H_k(X)\otimes {\mathbb{Q}}\rightarrow
H_k(X,{\mathbb{Q}})$ is an surjective, then $L_1H_k(X)\otimes
{\mathbb{Q}}\rightarrow H_k(X,{\mathbb{Q}})$ must be surjective.
Proposition 3.1 gives the needed surjectivity for $k\geq 5$ even if
$X$ is singular variety of dimension 3.

\qe

\noindent\textbf{The proof of Corollary 1.2:} By Corollary 1.1, we
only need to show that
$T_1H_4(X,{\mathbb{Q}})=G_1H_4(X,{\mathbb{Q}})$. By the assumption
and Poincar\'{e} duality, $H_4(X,{\mathbb{Q}})\cong
H_2(X,{\mathbb{Q}})\cong {\mathbb{Q}}$. Therefore,
$G_1H_4(X,{\mathbb{Q}})=H_4(X,{\mathbb{Q}})\cong {\mathbb{Q}}$ and
again by the commutative diagram
$$
\begin{array}{ccccc}
& L_2H_k(X)\otimes {\mathbb{Q}} &  \stackrel{s}{\rightarrow} & L_1H_k(X)\otimes {\mathbb{Q}}  &    \\

&  \downarrow \Phi_{2,k} &   &\downarrow \Phi_{1,k} &\\

& H_k(X,{\mathbb{Q}}) &\cong   & H_k(X,{\mathbb{Q}}), &

\end{array}
$$
we have the surjectivity of $L_1H_4(X)\otimes
{\mathbb{Q}}\rightarrow H_4(X,{\mathbb{Q}})$.

\qe

\noindent\textbf{The proof of Corollary 1.3:} Suppose $X=S\times C$,
where $S$ is a smooth projective surface and $C$ is a smooth
projective curve. We only need to consider  the surjectivity of
$L_1H_4(X)\otimes {\mathbb{Q}}\rightarrow H_4(X,{\mathbb{Q}})$
because of Corollary 1.1. Now the K\"unneth formula for the rational
homology of $H_4(S\times C, {\mathbb{Q}})$ and Theorem 2.1 for $S$
and $C$ gives the surjectivity in this case.

\qe

\noindent\textbf{The proof of Corollary 1.4:} This follows directly
from Theorem 1.3.

\qe

\noindent\textbf{The proof of Corollary 1.5:} By Theorem 1.4, we
only need to show that
$T_pH_k(X,{\mathbb{Q}})=G_pH_k(X,{\mathbb{Q}})$ for $k\geq 2p+2$. By
the definition of geometric definition, an element $a\in G_p
H_k(X,{\mathbb{Q}})$ comes from the linear combination of elements
$b_j\in H_k(Y_j,{\mathbb{Q}})$ for subvarieties $Y_j$ of dim$Y_j\leq
k-p$. From the following commutative diagram
$$
\begin{array}{ccccc}
 &i_*: L_pH_k(Y)\otimes {\mathbb{Q}}  &  \rightarrow & L_pH_k(X)\otimes {\mathbb{Q}}  &    \\

 &  \downarrow \Phi_{p,k} &   &\downarrow \Phi_{p,k} &\\

 & i_*: H_k(Y,{\mathbb{Q}}) &\rightarrow   & H_k(X,{\mathbb{Q}}), &

\end{array}
$$
it is enough to show that $L_pH_k(Y)\rightarrow H_k(Y)$ is
surjective for any irreducible subvariety $Y\subset X$ with dim$(Y)=
k-p$. By Suslin's conjecture, this is true for any smooth variety
$Y$ since $\dim(Y)=k-p$. Now we need to show that it is also true
for singular irreducible varieties if the Sulin Conjecture is true.

Using induction, we will show the following lemma:

{\Lemma If the Suslin Conjecture is true for every smooth projective
variety, then it is also true for every quasi-projective variety.}

\medskip
\bp Suppose that $Y$ is an irreducible quasi-projective
variety with $\dim(Y)=m$, $S$ is an irreducible quasi-projective
variety with $\dim(S)=n<m$ and

$$
\left\{
\begin{array}{l}
 L_pH_{n+p-1}(S)\rightarrow H_{n+p-1}(S) \quad is\quad injective,\\
 L_pH_{n+q}(S)\cong H_{n+q}(S) \quad for\quad q\geq p.
\end{array}
\right.
$$

Denote by $\overline{Y}$ a projective closure of $Y$ and
$S=sing(\overline{Y})$ the singular point set of $\overline{Y}$.
Let $U=\overline{Y}-S$
 Let $\sigma: \widetilde{Y}\rightarrow \overline{Y}$ be a desingularization of
$\overline{Y}$ and denote by $D:=\widetilde{Y}-U$. The existence of
a smooth $\widetilde{Y}$ is guaranteed  by Hironaka \cite{Hironaka}.
Then $D$ is the union of irreducible varieties with dimension $\leq
m-1$.

By Lemma 2.3, we have the following commutative diagram
$$
\begin{array}{ccccccccccc}
\cdots\rightarrow & L_pH_k(Z) & \rightarrow & L_pH_k({V})
& \rightarrow & L_pH_k(U) & \rightarrow & L_pH_{k-1}(Z) & \rightarrow & \cdots & \\

& \downarrow \Phi_{p,k}&   & \downarrow \Phi_{p,k}&   & \downarrow
\Phi_{p,k}&   & \downarrow \Phi_{p,k-1}&   &  &\\

 \cdots\rightarrow & H_k(Z) & \rightarrow &
H_k(V) & \rightarrow & H_k^{BM}(U)&\rightarrow & L_pH_{k-1}(Z) &
\rightarrow & \cdots,&

\end{array}
$$
where $U\subset V$ are quasi-projective varieties of
$\dim(V)=\dim(U)=m$ and $Z=V-U$ is a closed subvariety of $V$.

\medskip
\noindent\textbf{Claim:} By inductive assumption, the above
commutative diagram and the Five Lemma, we have the equivalence
between
$$
\left\{
\begin{array}{l}
 L_pH_{m+p-1}(U)\rightarrow H_{m+p-1}(U) \quad is\quad injective,\\
 L_pH_{m+q}(U)\cong H_{m+q}(U) \quad for\quad q\geq p.
\end{array}
\right.
$$
 and
$$
\left\{
\begin{array}{l}
 L_pH_{m+p-1}(V)\rightarrow H_{m+p-1}(V) \quad is\quad injective,\\
 L_pH_{m+q}(V)\cong H_{m+q}(V) \quad for\quad q\geq p.
\end{array}
\right.
$$
The proof of the claim is obvious.

By using the claim for finite times beginning from
$V=\widetilde{Y}$, we have the result for any quasi-projective
variety $U$. The proof of Lemma 3.2 is done.

\qe

By Lemma 3.2, we know that the Suslin's Conjecture is also true for
singular varieties. This completes the proof of Corollary 1.4.

\qe

\begin{center}{\bf Acknowledge}\end {center} I
would like to express my gratitude to my advisor, Blaine Lawson,
for all his help.

\medskip

\noindent
Department of Mathematics,\\
Stony Brook University, SUNY,\\
Stony Brook, NY 11794-3651\\
Email:wenchuan@math.sunysb.edu


\begin{thebibliography}{AAAA}
%
%
\bibitem[AKMW]{AKMW} Abramovich, Dan; Karu, Kalle; Matsuki, Kenji; W\l
odarczyk, Jaros\l aw, {\sl Torification and factorization of
birational maps.}  J. Amer. Math. Soc. 15 (2002), no. 3, 531--572
(electronic).

\bibitem[C]{Clemens} H. Clemens, {\sl Homological equivalence, modulo
algebraic equivalence, is not finitely generated.} Inst. Hautes
\'{E}tudes Sci. Publ. Math. No. 58 (1983), 19--38 (1984).

\bibitem[DT]{Dold-Thom} Dold, A. and Thom, R.,{\sl Quasifaserungen
und unendliche symmetrische Produkte.} (German) Ann. of Math. (2) 67
1958 239--281.


\bibitem[F1]{Friedlander1} Friedlander, Eric M., {\sl Algebraic cycles, Chow
varieties, and Lawson homology.}  Compositio Math. 77 (1991), no. 1,
55--93.

\bibitem[F2]{Friedlander2} Friedlander, Eric M.,{\sl Filtrations on algebraic
cycles and homology.} Ann. Sci. \'{E}cole Norm. Sup. (4) 28 (1995),
no. 3, 317--343.


\bibitem[FG]{Friedlander-Gabber} Eric M. Friedlander; Ofer Gabber,{\sl Cycle
spaces and intersection theory. Topological methods in modern
mathematics} (Stony Brook, NY, 1991), 325--370, Publish or Perish,
Houston, TX, 1993.

\bibitem[FHW]{Friedlander-Haesemesyer-Walker} Eric M. Friedlander,
Christian Haesemeyer, and Mark E. Walker, {\sl Techniques,
computations, and conjectures for semi-topological K-theory}
Preprint.

\bibitem[FL]{Friedlander-Lawson} Eric M. Friedlander; H.
Blaine Lawson, Jr. {\sl A theory of algebraic cocycles}. Ann. of
Math. (2) 136 (1992), no. 2, 361--428.

\bibitem[FM]{Friedlander-Mazur} Eric M. Friedlander; Barry Mazur,
{\sl Filtrations on the homology of algebraic varieties. With an
appendix by Daniel Quillen.}  Mem. Amer. Math. Soc. 110 (1994), no.
529, x+110 pp.


\bibitem[GH]{Griffiths-Harris} Griffiths, P.; Harris, J.,  {\sl Principles
of algebraic geometry.} Reprint of the 1978 original. Wiley Classics
Library.  John Wiley \& Sons, Inc., New York, 1994. xiv+813 pp. ISBN
0-471-05059-8

\bibitem[Gro]{Grothendieck} A. Grothendieck, {\sl Standard conjectures on algebraic
cycles}, Algebraic Geometry (Bombay, 1968), Oxford Univ. Press,
London, 1969, 193-199.

\bibitem[Hi]{Hironaka} H. Hironaka, {\sl Resolution of singularities
of an algebraic variety over a field of characteristic zero. I, II.}
Ann. of Math. (2) 79 (1964), 109--203; ibid. (2) 79 1964 205--326.

\bibitem[H]{author} W. Hu, {\sl Birational invariants defined by
Lawson homology.}  arXiv:math.AG/0511722.

\bibitem[L1]{Lawson1}
Lawson, H.B. Jr, {\sl Algebraic cycles and homotopy theory.}, Ann.
of Math. {\bf 129}(1989), 253-291.

\bibitem[L2]{Lawson2} Lawson, H. B. Jr, {\sl Spaces of algebraic
cycles.} pp. 137-213 in Surveys in Differential Geometry, 1995
vol.2, International Press, 1995.


\bibitem[Lew]{Lewis} Lewis, James D. {\sl A survey of the Hodge
conjecture. (English. English summary) Second edition. Appendix B by
B. Brent Gordon.} CRM Monograph Series, 10. American Mathematical
Society, Providence, RI, 1999. xvi+368 pp. ISBN 0-8218-0568-1


\bibitem[Lieb]{Lieberman} David I Lieberman, {\sl Numerical and homological equivalence of
algebraic cycles on Hodge manifolds.}  Amer. J. Math. 90 1968
366--374.


\bibitem[Li]{Lima-Filho} Lima-Filho, P., {\sl Lawson homology for
quasiprojective varieties.}  Compositio Math. 84(1992), no. 1,
1--23.


\bibitem[Wa]{Walker} Mark E. Walker, {\sl The morphic Abel-Jacobi map.}

\bibitem[W]{Wlodarczyk} W\l odarczyk, J., {\sl Toroidal varieties and
the weak factorization theorem.} Invent. Math. 154 (2003), no. 2,
223--331.
\end{thebibliography}
\end{document}